\input amstex
\documentstyle{amsppt}
\topmatter \magnification=\magstep1 \pagewidth{5.2 in}
\pageheight{6.7 in}
\abovedisplayskip=10pt \belowdisplayskip=10pt
\parskip=8pt
\parindent=5mm
\baselineskip=2pt
\title
 On the $q$-analogue of two-variable $p$-adic $L$-function
\endtitle
\author   Taekyun Kim    \endauthor

\affil{ {\it Institute of Science Education,\\
        Kongju National University, Kongju 314-701, S. Korea\\
        e-mail: tkim64$\@$hanmail.net ( or tkim$\@$kongju.ac.kr)}}\endaffil
\thanks  2000 Mathematics Subject Classification:  11S80, 11B68, 11M99 .\endthanks
\abstract{ We construct the two-variable $p$-adic $q$-$L$-function
which interpolates the generalized $q$-Bernoulli polynomials
associated with primitive Dirichlet character $\chi$. Indeed, this
function is the $q$-extension of two-variable $p$-adic
$L$-function due to Fox, corresponding to the case $q=1 .$
Finally, we give some $p$-adic integral representation for this
two-variable $p$-adic $q$-$L$-function and derive to $q$-extension
of the generalized formula of Diamond and Ferro and Greenberg for
the two-variable $p$-adic $L$-function in terms of the $p$-adic
gamma and $\log$ gamma function.
 }\endabstract
\rightheadtext{  Two-variable $p$-adic $q$-$L$-function}
\leftheadtext{T. Kim}
\endtopmatter

\document

\define\Z{\Bbb Z_p}

\head \S 1. Introduction \endhead Let $p$ be a fixed prime.
Throughout this paper  $\Bbb Z_p,\,\Bbb Q_p , \,\Bbb C$ and $\Bbb
C_p$ will, respectively, denote the ring of $p$-adic rational
integers, the field of $p$-adic rational numbers, the complex
number field and the completion of algebraic closure of $\Bbb Q_p
. $  Let $v_p$ be the normalized exponential valuation of $\Bbb
C_p$ with $|p|_p=p^{-v_p(p)}=p^{-1}.$ When one talks of
$q$-extension, $q$ is variously considered as an indeterminate, a
complex number $q\in\Bbb C ,$ or a $p$-adic number $q\in\Bbb C_p$.
If $q\in\Bbb C ,$ one normally assumes $|q|<1 .$  If $ q \in \Bbb
C_p ,$  then we assume $|q-1|_p < p^{-\frac1{p-1}},$ so that
$q^x=\exp(x\log q)$ for $|x|_p \leq 1.$  Kubota and Leopoldt
proved the existence of meromorphic functions, $L_p(s, \chi )$,
defined over the $p$-adic number field, that serve as $p$-adic
equivalents of the Dirichlet $L$-series, cf.[ 9, 11 ].
 These $p$-adic
$L$-functions interpolate the values
$$L_p(1-n, \chi)=-\frac{1}{n}(1-\chi_n(p)p^{n-1})B_{n,\chi_n},
\text{ for $n\in\Bbb N=\{1, 2,\cdots, \}  $ ,}$$ where
$B_{n,\chi}$ denote the $n$th generalized Bernoulli numbers
associated with the primitive Dirichlet character $\chi ,$ and
$\chi_n=\chi w^{-n} ,$ with $w$  the $Teichm\ddot{u}ller$
character, cf.[1-37]. In this paper, we use the notation:
$$[x]_q=[x:q]=\frac{1-q^x}{1-q}.$$
Hence,  $\lim_{q\rightarrow1}[x]=x$ for any $x$ with $|x|_p \leq 1
$ in the present $p$-adic case. Let $d$ be a fixed integer and let
$p$ be a fixed prime number. We set
$$\align
&X=X_d=\varprojlim_N (\Bbb Z/dp^N\Bbb Z),\text{ $X_1=\Bbb Z_p$ ,}\\
&X^*=\bigcup\Sb 0<a<dp\\ (a,p)=1\endSb a+dp\Bbb Z_p,\\
&a+dp^N\Bbb Z_p=\{x\in X\mid x\equiv a\pmod{dp^N}\},
\endalign$$
where $a\in \Bbb Z$ lies in $0\leq a<dp^N,$ cf. [11-30].

For any positive integer $N,$  we set $$\mu_q(a+dp^N
\Z)=\frac{q^a}{[dp^N]} , \text{ cf.[18] } ,$$  and this can be
extended to a distribution on $X.$ This distribution yields an
integral for each non-negative integer $m$:
$$ \beta_{m,q}^{*}=\int_{\Bbb Z_p} [x]^m d\mu_q(x)=\int_{X}[a]^{m}d\mu_q(a)=\frac1{(1-q)^m}
\sum_{i=0}^m\binom mi(-1)^i \frac{i+1}{[i+1]} ,$$ where
$\beta_{m,q}^{*}$ are the $m$th Carlitz's $q$-Bernoulli numbers,
cf. [3, 4, 5 ]. In a recent paper [ 10], Fox defined a
two-variable $p$-adic $L$-function $L_p(s, t|\chi)$ with the
property that
$$L_p(1-m, t|\chi)=-\frac{B_{m,\chi_m}(p^*
t)-\chi_m(p)p^{m-1}B_{m,\chi_m}(p^{-1}p^* t)}{m}, $$ for positive
integer $m$ and $t\in\Bbb C_p$ with $|t|_p\leq 1 ,$ where $p^*=p$
if $p>2$ and $p^*=4$ if $p=2 ,$ and $B_{m,\chi_m}(x)$ are the
$m$th generalized Bernoulli polynomials attached to $\chi .$ In
[12 ], we gave the interesting results on the $p$-adic
$q$-$L$-functions, a subject initiated by Neal Koblitz [ 24], in
the beginning of the 1980's in which the author made some
contributions, cf.[11-23 ]. For positive integer $n ,$ these
functions satisfy $$L_{p,q}(1-n, \chi)=-\frac{\beta_{n,q,
\chi_n}-[p]_q^{n-1}\chi_n(p)\beta_{n,q^p,\chi_n}}{n},$$ where
$\beta_{n, q,\chi_n}$ are the $n$th generalized $q$-Bernoulli
numbers attached to $\chi$ which are defined by author, cf.[ 20].
The purpose of this paper is to construct a two-variable $p$-adic
$q$-$L$-function $L_{p,q}(s, t|\chi)$ for the Dirichlet character
$\chi$ with the property that
$$L_{p,q}(1-m, t|\chi)=-\frac{\beta_{n,q,\chi_n}(p^*
t)-\chi_{n}(p)[p]_q^{n-1}\beta_{n,q^p,\chi_n}(p^{-1}p^* t)}{n}
,\text{ $n\in\Bbb Z_{+}$,}$$ where $\beta_{n,q,\chi_n}(x)$ are the
$n$th generalized $q$-Bernoulli polynomials attached to $\chi$.
This function is actually a $q$-extension of the two-variable
$p$-adic $L$-function of Fox, corresponding to the case $q=1.$ For
a prime number $p$ and for a Dirichlet character defined modulo
some integer, the $p$-adic $L$-function was constructed by
interpolating the values of complex analytic $L$-function at
non-positive integers. Diamond [7, 8 ] obtained formulas which
express the values of $p$-adic $L$-function at positive integers
in terms of the $p$-adic $\log$ gamma function. In this paper, we
give the $q$-extension of his results to the case of the
two-variable $p$-adic $q$-$L$-function and obtain the formulas
 which express the values of $\frac{\partial}{\partial s}L_{p,q}(0,t|\chi)$
 in terms of the $q$-extension of Diamond $p$-adic $\log$ gamma function. Finally,
we give the values of  $L_{p,q}(s,t|\chi)$ at $s=1$.

 \head 2. $q$-extension of two-variable Dirichlet's $L$-series  \endhead

In this section we assume that $q\in\Bbb C$ with $|q|<1 .$ The
$q$-Bernoulli numbers are usually defined  by
$$\beta_{0,q}=\frac{q-1}{\log q}, \quad (q\beta_q +1)^n
-\beta_{n,q}=\delta_{n,1}, \tag1$$ where $\delta_{n,1}$ is the
Kronecker symbol and we use the usual convention about replacing
$\beta_q^i$ by $\beta_{i,q} $, cf.[13, 14 ]. Note that
$\lim_{q\rightarrow 1}\beta_{k,q}=B_k ,$ where $B_k$ are the $k$th
ordinary Bernoulli numbers. In [14] the $q$-Bernoulli polynomials
are defined by
 $$\beta_{n,q}(x)=\sum_{i=0}^n \binom ni
 q^{xi}\beta_{i,q}[x]_q^{n-i}=\frac{1}{(1-q)^n}\sum_{i=0}^n\binom
 ni \frac{i}{[i]_q}(-q^x)^i. \tag2$$
 From the Eq.(1), we note that
 $$\beta_{n,q}=\frac{1}{(1-q)^n}\sum_{i=0}^n\binom
 ni(-1)^{n-i}\frac{i}{[i]_q}, \text{ where $\binom ni$ is the binomial
 coefficient.}$$
Thus, we have the generating function of $q$-Bernoulli numbers as
follows:
$$F_q(t)=e^{\frac{t}{1-q}}\sum_{j=0}^{\infty}\frac{j}{[j]_q}\left(\frac{1}{q-1}\right)^j\frac{t^j}{j!}
=\sum_{j=0}^{\infty}\frac{\beta_{j,q}}{j!}t^j, \text{ for $|t|<1$
}. \tag 3$$ By (3), we easily see that the $q$-Bernoulli numbers
are the unique solutions of the following $q$-difference equation
in the complex plane:
$$F_q(t)=\frac{q-1}{\log
q}e^{\frac{t}{1-q}}-t\sum_{n=0}^{\infty}q^ne^{[n]_qt}, \text{ for
$|t|<1$}. \tag4$$ In the Eq.(2), we consider the generating
function of $q$-Bernoulli polynomials as follows:
$$\sum_{n=0}^{\infty}\beta_{n,q}(x)\frac{t^n}{n!}=F_q(x,t)
=e^{\frac{t}{1-q}}\sum_{j=0}^{\infty}\frac{j}{[j]_q}\left(\frac{1}{q-1}\right)^j
q^{jx}\frac{t^j}{j!}, \text{ for $|t|<1$. }\tag 5$$ Thus, we
obtain the following $q$-difference equation for the generating
function of $q$-Bernoulli polynomials in the complex plane:
$$F_q(x,t)=\frac{q-1}{\log
q}e^{\frac{t}{1-q}}-t\sum_{n=0}^{\infty}q^{n+x}e^{[n+x]_q t},
\text{ for $|t|<1$. } \tag6 $$ Let $\chi$ be the Dirichlet
character with conductor $f=f_{\chi}\in\Bbb N.$ Then the
generalized $q$-Bernoulli numbers attached to $\chi$,
$\beta_{n,\chi, q} ,$ are defined by
$$F_{q,\chi}(t)=-t\sum_{n=1}^{\infty}\chi(n)q^ne^{[n]_q
t}=\sum_{n=0}^{\infty}\beta_{n,\chi,q}\frac{t^n}{n!}, \quad |t|<1.
\tag7$$ Remark. From the Eq.(7), we note that
$$\lim_{q\rightarrow 1}F_{q,\chi}(t)=\sum_{a=1}^f \frac{\chi(a)e^{at}
t}{e^{ft}-1}=\sum_{n=0}^{\infty}B_{n,\chi}\frac{t^n}{n!}.$$ By the
Eq.(7), we easily see that
$$\beta_{n,\chi,q}=[f]_q^{k-1}\sum_{a=1}^f\chi(a)\beta_{k,q^f}(\frac{a}{f}).
\tag8$$ We now also define the generalized $q$-Bernoulli
polynomials attached to $\chi$ as follows:
$$F_{q,\chi}(x,t)=-t\sum_{n=1}^{\infty}\chi(n)q^{n+x}e^{[n+x]_q
t}=\sum_{n=0}^{\infty}\beta_{n,\chi,q}(x)\frac{t^n}{n!}.\tag9$$
Thus, we obtain the below formula:
$$F_{q,\chi}(x,t)=-t\sum_{a=1}^f
\chi(a)\sum_{n=0}^{\infty}q^{nf+a+x}e^{[nf+a+x]_q t}. $$ From
this, we note that
$$\beta_{n,\chi,q}(x)=\sum_{k=0}^n\binom nk
q^{kx}\beta_{k,\chi}[x]_q^{n-k}=(q^x\beta_{\chi}+[x]_q)^n ,$$ with
the usual convention about replacing $\beta_{\chi}^n$ by
$\beta_{n,\chi,q}.$

Let $g$ be a positive integral multiple of $f=f_{\chi}$. Then for
each $n\in\Bbb Z$, $n\geq 0 ,$ we have
$$\aligned
&\sum_{n=0}^{\infty}\beta_{n,\chi,q}(x)\frac{t^n}{n!}\\
&=-t\sum_{a=0}^{g-1}\chi(a)\sum_{n=0}^{\infty}q^{gn+a+x}e^{[gn+a+x]_q
t}=-t\sum_{a=0}^{g-1}\chi(-a)\sum_{n=1}^{\infty}q^{g(n+\frac{x-a}{g})}e^{[g]_q[\frac{x-a}{g}]_{q^g}t}\\
&=\sum_{n=0}^{\infty}\left([g]_q^{n-1}\sum_{a=0}^{g-1}\chi(-a)\beta_{n,q^g}(\frac{x-a}{g})\right)\frac{t^n}{n!}.
\endaligned$$
Therefore, we obtain the following lemma:
 \proclaim{ Lemma 1}
Let $g$ be a positive integral multiple of $f=f_{\chi}$. Then for
each $n\in\Bbb N \bigcup \{0\},$ we have
$$\beta_{n,\chi,q}(x)=[g]_q^{n-1}\sum_{a=0}^{g-1}\chi(a)\beta_{n,q^g}(\frac{x+a}{g})
=[g]_q^{n-1}\sum_{a=0}^{g-1}\chi(-a)\beta_{n,q^g}(\frac{x-a}{g}).$$\endproclaim
Note that the series on the right hand side of (6) and (9) are
uniformly convergent. Hence, we easily see that
$$\beta_{k,q}(x)=\frac{d^k}{dt^k}F_q(x,t)|_{t=0}
=\frac{q-1}{\log q}\frac{1}{(1-q)^k}-k
\sum_{n=0}^{\infty}q^{n+x}[n+x]_q^{k-1}, $$ and,
$$ \beta_{k,\chi,q}(x)=\frac{d^k}{dt^k}F_{q,\chi}(x,t)|_{t=0}
=-k\sum_{n=1}^{\infty}[n+x]_q^{k-1}q^{n+x}\chi(n), \quad k\geq
1.$$
 Therefore, we obtain the following theorem:
\proclaim{Theorem 2} Let $\chi $ be the primitive character with
conductor $f=f_{\chi}$. For $k\geq 1 ,$ $q\in\Bbb C$ with $|q|<1
,$ we have
$$\beta_{k,q}(x)=\frac{q-1}{\log q}\frac{1}{(1-q)^k}-k
\sum_{n=0}^{\infty}q^{n+x}[n+x]_q^{k-1},$$ and,
$$\beta_{k,\chi,q}(x)=-k\sum_{n=1}^{\infty}[n+x]_q^{k-1}q^{n+x}\chi(n).$$
\endproclaim

In [ ], the $q$-analogue of the Hurwitz's zeta function was
defined by
$$\zeta_q(s,x)=\sum_{n=0}^{\infty}\frac{q^{n+x}}{[n+x]_q^s}-\frac{1}{(s-1)}\frac{(1-q)^s}{\log
q}, \quad s\in\Bbb C. \tag10$$ This function is meromorphic for
$s\in\Bbb C$ with simple pole at $s=1 .$ By using Theorem 2, we
easily see that $\zeta_q(1-n, x)=-\frac{\beta_{n,q}(x)}{n},$
$n\in\Bbb N$.
 From the results of the above Theorem 2, we  can consider the
 $q$-analogue of two-variable Dirichlet's $L$-series as follows:
\proclaim{ Definition 3} For $s\in\Bbb C ,$ we define the
$q$-analogue of two-variable Dirichlet $L$-series as
$$L_q(s,x|\chi)=\sum_{n=1}^{\infty}q^{n+x}\frac{\chi(n)}{[n+x]_q^s}
=(q-1)\sum_{n=1}^{\infty}\frac{\chi(n)}{[n+x]_q^{s-1}}+\sum_{n=1}^{\infty}\frac{\chi(n)}{[n+x]_q^s}.$$
\endproclaim
In [12] one variable $q$-$L$-series are defined by
$L_q(s,\chi)=\sum_{n=1}^{\infty}\frac{\chi(n)}{[n]_q^s} ,$
$s\in\Bbb C .$ Thus, we see that $L_q(s,0|\chi)=L_q(s,\chi) .$ In
the previous paper [23] we also defined the two-variable Dirichlet
$L$-series  as
$L(s,x|\chi)=\sum_{n=1}^{\infty}\frac{\chi(n)}{(n+x)^s}, $
$s\in\Bbb C. $ From this definition, we note that
$\lim_{q\rightarrow 1} L_q(s,x|\chi)=L(s,x|\chi).$ By Theorem 2
and Definition 3, we obtain the following corollary:

 \proclaim{Corollary 4}
 If $k\geq 1,$ $0\leq x\leq 1$, then we have
 $$ L_q(1-k,x|\chi)=-\frac{\beta_{k,\chi, q}(x)}{k}.$$
\endproclaim
Remark. Let $\chi$ be a Dirichlet character with conductor
$f=f_{\chi} .$ By (7) and (8), we easily see that
$$\sum_{k=0}^{nf-1}\chi(k)q^k[k]_q^l=\frac{1}{l+1}\left(
\beta_{l+1,\chi,q}(nf)-\beta_{l+1,\chi,q}\right), \quad n,
l\in\Bbb N . \tag11$$ In [31] M. Schlosser investigated the
$q$-analogues of the sums of consecutive integers, squares, cubes,
quarts and quints. That is, his q-analogue of $\sum_{k=1}^n k^m ,$
for $m=1, 2, 3, 4, 5 ,$ gave by employing specific identities for
very-well-posed basic hypergeometric series, in conjunction with
suitable specialization of the parameters. In the final page of
his paper he did guess that reasonable continuation involving
higher integer powers will follow the same his pattern. The
Eq.(11) is the generalization for the problem which was guessed by
Schlosser in [31]. Indeed if we take trivial character in Eq.(11),
then Eq.(11) becomes the $q$-analogue of the sums of powers of
consecutive integers involving higher order. By using the Eq.(10)
and Definition 3, we obtain the below identity:
$$L_q(s,x|\chi)=[f]_q^{-s}\sum_{a=1}^f \chi(a)\zeta_{q^f}(s,
\frac{a+x}{f}). \tag12$$ Let $\Gamma(s)$ be the gamma function.
Then we can readily see that

$$\aligned
&\frac{1}{\Gamma(s)}\int_0^{\infty}t^{s-2}F_q(x,-t)dt\\
&=\frac{q-1}{\log q}\frac{1}{\Gamma(s)}
\int_{0}^{\infty}(1-q)^{s-1}x^{s-1}e^{-x}dx+\sum_{n=0}^{\infty}
\frac{q^{n+x}}{[n+x]_q^s}\frac{1}{\Gamma(s)}\int_0^{\infty}e^{-xt}t^{s-1}dt\\
&=-\frac{1}{s-1}\frac{(1-q)^s}{\log q} +
\sum_{n=0}^{\infty}\frac{q^{n+x}}{[n+x]_q^s}=\zeta_q(s,x).\endaligned$$
Therefore, we obtain the followings lemma:

\proclaim{ Lemma 5} For $s \in\Bbb C ,$ we have
$$\zeta_q(s,x)=\frac{1}{\Gamma(s)}\int_{0}^{\infty}t^{s-2}\left(
\frac{q-1}{\log q} e^{-\frac{t}{1-q}}+t\sum_{n=0}^{\infty}q^{n+x}
e^{-t[n+x]_q} \right)dt.$$
\endproclaim
By applying Mellin transforms  and residue theorem to Lemma 5, we
can derive the below identity:
$$\zeta_q(1-n, x)=\frac{(-1)^n}{n!}\beta_{n,q}(x)(2\pi
i)\left(\frac{(n-1)!}{2\pi i}(-1)^{n-1}\right), \quad n\in\Bbb N.
$$
In the case $x=1$ we note that
$\zeta_q(1-n,1)=-\frac{\beta_{n,q}(1)}{n}=-\frac{(q\beta_q
+1)^n}{n}=-\frac{\beta_{n,q}}{n}, \quad n> 1 . $ By (9), we easily
see that $$\aligned
&\frac{1}{\Gamma(s)}\int_0^{\infty}F_{q,\chi}(x, -t)t^{s-2}dt\\
&=\sum_{n=1}^{\infty}\chi(n)q^{n+x}\frac{1}{\Gamma(s)}\int_{0}^{\infty}
e^{-[n+x]_q t}t^{s-1}dt=L_q(s,x|\chi), \quad s\in\Bbb
C.\endaligned \tag13$$ Thus, we have
$$L_q(1-n, x|\chi)=-\frac{\beta_{n,\chi,q}(x)}{n}, \quad n\in\Bbb N.$$
Let
$$\aligned
&H_q(s,a,F)=\sum_{\Sb m\equiv a(\mod F)\\m>0 \endSb}
\frac{q^m}{[m]_q^s}+\frac{1}{F}\frac{(1-q)^s}{(1-s)\log q}\\
&=\sum_{n=0}^{\infty}\frac{q^{a+nF}}{[a+nF]_q^s}+\frac{1}{F}\frac{(1-q)^s}{(1-s)\log
q}=[F]_q^{-s}\zeta_{q^F}(s,\frac{a}{F}),
 \endaligned\tag14$$
 where $a$ and $F$ are positive integers with $0<a<F. $
Let $\chi (\neq 1)$ be the Dirichlet character with conductor $F$.
Then the $q$-analogue of Dirichlet $L$-function can be expressed
as the sum
$$L_q(s,\chi)=\sum_{a=1}^{F}\chi(a)H_q(s,a,F), \quad s\in\Bbb C,
\text{ cf.[17].}\tag15$$ The function $H_q(s, a, F)$ is a
meromorphic for $s\in\Bbb C$ with simple pole at $s=1$, having
residue $\frac{q-1}{F\log q}, $ and it interpolates the values
$$H_q(1-n, a, F)=-\frac{[F]_q^{n-1}}{n}\beta_{n,q^F}(\frac{a}{F}),
\quad\text{where $n\in\Bbb Z,$ $n\geq 1$}.\tag16$$ We now modify
the partial $q$-zeta function as follows:
$$H_q(s,a,
F)=\frac{1}{s-1}\frac{1}{[F]_q}[a]_q^{1-s}\sum_{j=0}^{\infty}\binom{1-s}jq^{aj}\beta_{j,q^F}
\left[\frac{F}{a}\right]_{q^a}^j ,\quad \text{ for $s\in\Bbb C$}.
\tag 17$$ By (14), (15), (16) and (17), we easily see that
$$L_q(s,\chi)=\frac{1}{s-1}\frac{1}{[F]_q}\sum_{a=1}^{F}\chi(a)[a]_q^{1-s}
\sum_{m=0}^{\infty}\binom{1-s}m
q^{am}\beta_{m,q^F}\left[\frac{F}{a}\right]_{q^a}^m. \tag18$$ From
the Definition 3, we note that $L_q(s, x|\chi)$ is analytic for
$s\in\Bbb C ,$ except $s\neq 1$ when $\chi \neq 1 .$ Using Eq.(17)
to define $H_q(s,a+x, F)$ for all $a \in \Bbb Z$ with $0<a<F$,
$x\in\Bbb R$ with $0\leq x\leq 1 ,$ we obtain
$$L_q(s,x|\chi)=\sum_{a=1}^{F}\chi(a)H_q(s,a+x, F). \tag19 $$
Let $F$ and $a$ be positive integers with $0<a<F ,$ and let
$$L_q(s,x|\chi)=\frac{1}{s-1}\frac{1}{[F]_q}\sum_{a=1}^F\chi(a)[a+x]_q^{1-s}
\sum_{m=0}^{\infty}\binom{1-s}{m}q^{(a+x)m}\beta_{m,q^F}\left[\frac{F}{a+x}\right]_{q^{a+x}}^m.\tag20$$
Then, $L_q(s,x|\chi)$ is analytic for $s\in\Bbb C ,$ except $s\neq
1$ when $\chi =1$. Furthermore, for each $n\in\Bbb Z$ with $n\geq
1$, we have
$$L_q(1-n, x|\chi)=-\frac{\beta_{n,\chi, q}(x)}{n}.$$ In this
section we introduced some of the basic facts about one-variable
$q$-$L$-series and two-variable $q$-$L$-series in complex plane.
Then their values at negative integers are given in terms of
generalized $q$-Bernoulli numbers and polynomials attached to
$\chi$. We also evaluate $L_q(1,x|\chi)$ and give some relation
with $q$-Bernoulli numbers and polynomials. By the definition of
$L_q(s,x|\chi)$, we easily see that
$$\aligned
&L_q(s,x|\chi)\\
 &=\frac{1}{s-1}\frac{1}{[F]_q}\sum_{a=1}^F\chi(a)
([a+x]_q^{1-s}+[a+x]_q^{1-s}\sum_{m=1}^{\infty}\binom{1-s}{m}
q^{(a+x)m}\beta_{m,q^F}[\frac{F}{a+x}]_{q^{a+x}}^m ).\endaligned$$
We now give the below Taylor expansion of $[a+x]_q^{1-s}$ at
$s=1$: $$[a+x]_q^{1-s}=1-(s-1)\log [a+x]_q+\cdots .$$ Thus, we see
that
$$\aligned
&L_q(1,x|\chi)\\
&=\frac{1}{[F]_q}\sum_{a=1}^F\chi(a)\left\{-\log ([a+x]_q)
+\sum_{m=1}^{\infty}\frac{(-1)^m}{m}q^{(a+x)m}\beta_{m,q^F}\left[\frac{F}{a+x}\right]_{q^{a+x}}^m\right\}.
\endaligned$$
In the case $x=0$, we have
$$
L_q(1,0|\chi)=L_q(1,\chi)
=\frac{1}{[F]_q}\sum_{a=1}^F\chi(a)\left\{-\log ([a]_q)
+\sum_{m=1}^{\infty}\frac{(-1)^m}{m}q^{am}\beta_{m,q^F}\left[\frac{F}{a}\right]_{q^{a}}^m\right\}.$$
The values of $L_q(s,x|\chi)$ at negative integers are algebraic,
hence may be regarded as lying in an extension of $\Bbb Q_p$. We
therefore look for a $p$-adic function which agrees with
$L_q(s,x|\chi)$ at the negative integers in the next section.

 \head 3. Two-variable $p$-adic $q$-$L$-functions  \endhead

In this section we shall consider the $p$-adic  analogs of the
two-variable $q$-$L$-functions which were introduced in the
previous section. Indeed this functions are the $q$-analogs of the
$p$-adic functions due to Fox, corresponding to the case $q=1$.
Let $w$ denote  the $Teichm\ddot{u}ller$ character, having
conductor $f_w =p^*$.  For an arbitrary character $\chi$, we
define $\chi_n=\chi w^{-n},$ where $n\in\Bbb Z$, in the sense of
the product of characters. Throughout this section, we assume that
$q\in\Bbb C_p$ with $|1-q|_p <p^{-\frac{1}{p-1}} .$ Let
$<a>=<a:q>=w^{-1}(a)[a]_q=\frac{[a]_q}{w(a)}.$ Then, we note that
$<a>\equiv 1$ $(\mod p^*p^{\frac{1}{p-1}}).$ By the definition of
$<a>$, we easily see that
$<a+p^*t>=w^{-1}(a+p^*t)[a+p^*t]_q=w^{-1}(a)[a]_q+w^{-1}(a)q^a[p^*t]_{q^a}\equiv
1$ $(\mod p^*p^{\frac{1}{p-1}}),$ where $t\in\Bbb C_p$ with
$|t|_p\leq 1.$ The $p$-adic logarithm function, $\log_p ,$ is the
unique function $\Bbb C_p^{\times} \rightarrow \Bbb C_p $ that
satisfy  (1) $\log_p (1+x)=\sum_{n=1}^{\infty}\frac{(-1)^n}{n}x^n
,$ $|x|_p<1$, (2) $\log_p(xy)=\log_p x+\log_p y ,$ $\forall x,
y\in\Bbb C_p^{\times},$ and $\log_p p =0 .$ Let
$A_j(x)=\sum_{n=0}^{\infty}a_{n,j}x^n$, $a_{n,j}\in\Bbb C$, $j=0,
1, 2,\cdots$ be a sequence of power series, each of which
converges in a fixed subset $D=\{s\in\Bbb C_p||s|_p\leq
|p^*|^{-1}p^{-\frac{1}{p-1}}\}$ of $\Bbb C_p$ such that (1)
$a_{n,j}\rightarrow a_{n, 0}$ as $j\rightarrow \infty$ for
$\forall n$; (2) for each $s\in D$ and $\epsilon >0$, there exists
$n_0=n_0(s,\epsilon)$ such that $\left|\sum_{n\geq
n_0}a_{n,j}s^n\right|_p<\epsilon $ for $\forall j$. Then
$\lim_{j\rightarrow \infty}A_j(s)=A_0(s)$ for all $s\in D .$ This
is used by Washington [36] to show that each of the function
$w^{-s}(a)a^s$ and $\sum_{m=0}^{\infty}\binom sm
\left(\frac{F}{a}\right)^mB_m ,$ where $F$ is the multiple of
$p^*$ and $f=f_{\chi}$, is analytic in $D$. Let $F$ be a positive
integral multiple of $p^*$ and $f=f_{\chi}$, and let
$$\aligned
&L_{p,q}(s,t|\chi)\\
&=\frac{1}{s-1}\frac{1}{[F]_q}\sum_{\Sb a=1 \\
(a,p)=1\endSb}^F\chi(a)<a+p^*t>^{1-s}
\sum_{m=0}^{\infty}\binom{1-s}m
\beta_{m,q^F}q^{(a+p^*t)m}[\frac{F}{a+p^*t}]_{q^{a+p^*t}}^m.\endaligned$$
Then $L_{p,q}(s,t|\chi)$ is analytic for $t\in\Bbb C_p$ with
$|t|_p\leq 1 ,$ provided $s\in D$, except $s\neq 1$ when $\chi\neq
1.$ For $t\in\Bbb C_p$ with $|t|_p\leq 1$, we see that
$\sum_{j=0}^{\infty}\binom sj
\beta_{j,q^F}q^{(a+p^*t)j}\left[\frac{F}{a+p^*t}\right]_{q^{a+p^*t}}^j$
is analytic for $s\in D .$ It readily follows that
$<a+p^*t>^s=w^{-s}(a)[a+p^*t]_q^s=<a>^s\sum_{m=0}^{\infty}\binom
sm \left(q^a[a]_q^{-1}[p^*t]_q\right)^m $ is analytic for
$t\in\Bbb C_p$ with $|t|_p\leq 1$ when $s\in D .$ Thus, since
$(s-1)L_{p,q}(s,t|\chi)$ is a finite sum of products of these two
functions, it must also be analytic for $t\in\Bbb C_p$, $|t|_p\leq
1$, whenever $s\in D.$ Note that
$$\lim_{s\rightarrow 1}(s-1)L_{p,q}(s,t|\chi)=\cases
\frac{1}{[F]_q}\frac{q^F-1}{\log q}(1-\frac{1}{p}),&\text{ if }
\chi =1,\\
0,&\text{ if } \chi \neq 1.\endcases$$

We now let $n\in\Bbb Z,$ $ n\geq 1, $ and fix $t\in\Bbb C_p$ with
$|t|_p\leq 1$. Since $F$ must be a multiple of $f=f_{\chi_n}$,
Lemma 1 implies that
$$\beta_{n,\chi_n ,
q}(p^*t)=[F]_q^{n-1}\sum_{a=0}^{F-1}\chi_n(a)\beta_{n,q^F}(\frac{a+p^*t}{F}).
\tag 21$$ If $\chi_n(p)=0$, then $(p, f_{\chi_n})=1$, so that
$\frac{F}{p}$ is a multiple of $f_{\chi_n}$. Therefore, we obtain
$$\chi_n(p)[p]_q^{n-1}\beta_{n,\chi_n, q^p}(p^{-1}p^*t)
=[F]_q^{n-1}\sum_{\Sb a=0\\p|a
\endSb}^F\chi_n(a)\beta_{n,q^F}(\frac{a+p^*t}{F}). \tag 22$$
The difference of these quantities  yields
$$\beta_{n,\chi_n,q}(p^*t)-\chi_n(p)[p]_q^{n-1}\beta_{n,\chi_n,q^p}(p^{-1}p^*t)
=[F]_q^{n-1}\sum_{\Sb a=1\\ p\nmid
a\endSb}^F\chi_n(a)\beta_{n,q^F}(\frac{a+p^*t}{F}).$$ By using the
distribution of $q$-Bernoulli polynomials,  we easily see that
$$\beta_{n,q^F}(\frac{a+p^*t}{F})=[F]_q^{-n}[a+p^*t]_q^n\sum_{m=0}^n
\binom nm q^{(a+p^*t)m}
\left[\frac{F}{a+p^*t}\right]_{q^{a+p^*t}}^m\beta_{m, q^F}. $$
Since $\chi_n(a)=\chi(a)w^{-n}(a)$ and for $(a,p)=1$, and
$t\in\Bbb C_p $ with $|t|_p\leq 1$, we have
$$\aligned
&\beta_{n,\chi_n,q}(p^*t)-\chi_n(p)[p]_q^{n-1}\beta_{n,\chi_n,q^p}(p^{-1}p^*t)\\
&=\frac{1}{[F]_q}\sum_{\Sb a=1\\ (a,p)=1\endSb
}^F\chi(a)<a+p^*t>^n\sum_{m=0}^{\infty}\binom nm
q^{(a+p^*t)m}\left[\frac{F}{a+p^*t}\right]_{q^{a+p^*t}}^m\beta_{m,q^F}
\endaligned$$
Thus, we see that
$$-\frac{1}{n}\left(\beta_{n,\chi_n,q}(p^*t)-\chi_n(p)[p]_q^{n-1}\beta_{n,\chi_n,q^p}(p^{-1}p^*t)\right)
=L_{p,q}(1-n,t|\chi),\text{ for $ n\in\Bbb N.$} $$ Therefore we
obtain the following theorem:

\proclaim{ Theorem 6} Let $F$ be a positive integral multiple of
$p^*$ and $f=f_{\chi}$, and let
$$\aligned
&L_{p,q}(s,t|\chi)\\
&=\frac{1}{s-1}\frac{1}{[F]_q}\sum_{\Sb
a=1\\(a,p)=1\endSb}^F\chi(a)<a+p^*t>^{1-s}\sum_{m=0}^{\infty}\binom{1-s}m
q^{(a+p^*t)m}\beta_{m,q^F}[\frac{F}{a+p^*t}]_{q^{a+p^*t}}^m
.\endaligned \tag 23$$ Then, $L_{p,q}(s, t|\chi)$ is analytic for
$t\in\Bbb C_p$, $|t|_p\leq 1,$ provided $s\in D$, except $s\neq 1$
when $\chi\neq 1 .$ Also, if $t\in\Bbb C_p ,$ $|t|_p\leq 1$, this
function is analytic for $s\in D$ when $\chi \neq 1$, and
meromorphic for $s\in D$, with simple pole at $s=1$ having residue
$\frac{1}{[F]_q}\frac{q^F-1}{\log q}(1-\frac{1}{p}) $ when
$\chi=1$.
Furthermore, for each $n\in\Bbb Z$, $n\geq 1$, we have
$$L_{p,q}(1-n,t|\chi)=-\frac{1}{n}\left(\beta_{n,\chi_n,q}(p^*t)-\chi_n(p)[p]_q^{n-1}
\beta_{n,\chi_n,q^p}(p^{-1}p^*t)\right).$$
\endproclaim
\proclaim {Remark} (1) Note that $L_{p,q}(s,
0|\chi)=L_{p,q}(s,\chi)$ for $s\in D$ with $s\neq 1$ if $\chi =1,$
where $L_{p,q}(s,\chi)$ is $p$-adic $q$-$L$-function, cf.[12].

(2) Let $L_p(s,t|\chi)$ be the two-variable $p$-adic $L$-functions
of Fox. Then we see that $\lim_{q\rightarrow 1}L_{p,q}(s,
t|\chi)=L_p(s,t|\chi) .$
\endproclaim

By means of a method provided by Washington [36], we now
generalize to two-variable $p$-adic $q$-$L$-function, $L_{p,q}(s,
t|\chi) ,$ by modifying $L_{p,q}(s,\chi),$ which was first defined
by the function
$$H_{p,q}(s,a,F)=\frac{1}{s-1}\frac{1}{[F]_q}<a>^{1-s}\sum_{j=0}^{\infty}
\binom{1-s}j \beta_{j,q^F}q^{aj}\left[\frac{F}{a}\right]_{q^a}^j,
$$ where $s\in D$, $s\neq 1$, $a\in\Bbb Z$ with $(a,p)=1$, and $F$
is a multiple of $p^*$, cf. [17].

The function $L_{p,q}(s,\chi)$ can be rewritten as the sum
$$L_{p,q}(s,\chi)=\sum_{\Sb a=1\\(a,p)=1
\endSb}^F\chi(a)H_{p,q}(s,a,F), \text{ cf.[17], }\tag24$$
provided $F$ is a multiple of both $p^*$ and $f=f_{\chi} .$ The
function $H_{p,q}(s,a, F)$ is a meromorphic for $s\in D$ with a
simple pole at $s=1$, having residue
$\frac{1}{[F]_q}\frac{1}{F}\frac{q^F-1}{\log q} ,$ and it
interpolates the values
$$H_{p,q}(1-n,a,F)=-\frac{1}{n}w^{-n}(a)\beta_{n,q^F}(\frac{a}{F}),
$$ where $n\in\Bbb Z$, $n\geq 1$, cf. [17, 12].

By using $H_{p,q}(s, a+p^*t, F)$, we can express
$L_{p,q}(s,t|\chi)$ for all $a\in\Bbb Z$, $(a,p)=1$, and $t\in\Bbb
C_p$ with $|t|_p\leq 1$, as follows:
$$L_{p,q}(s,t|\chi)=\sum_{\Sb a=1\\(a,p)=1
\endSb}^F\chi(a)H_{p,q}(s,a+p^*t, F). \tag 25$$
From the proof of Theorem 6, we note that $H_{p,q}(s, a+p^*t, F)$
is analytic for $t\in\Bbb C_p$, $|t|_p\leq 1$, where $s\in D$,
$s\neq 1$, and meromorphic for $s\in D$, with a simple pole at
$s=1$, when $t\in\Bbb C_p$, $|t|_p\leq 1$. Let us consider the
first partial derivative of the function $L_{p,q}(s,t|\chi)$ at
$s=0$. It is easy to see that
$$\frac{\partial^n}{\partial t^n }L_{p,q}(s,t|\chi)=\binom{-s}{n}n!
\left(p^*\frac{\log q}{q-1}\right)^nL_{p,q}(s+n,t|\chi_n),$$ for
all $s\in D$, $s\neq 1$ if $\chi =1$, and $t\in\Bbb C_p$ with
$|t|_p\leq 1 .$

Furthermore, we note that
$$\aligned
&\lim_{s\rightarrow 1-n}\binom{-s}n
L_{p,q}(s+n,t|\chi_n)=-\frac{(n-1)!}{n!}\lim_{s\rightarrow
1-n}(s+n-1)L_{p,q}(s+n,t|\chi_n)\\
&=-\frac{1}{n}\frac{1}{[F]_q}\sum_{\Sb a=1 \\(a,p)=1
\endSb}^F\chi_n(a)\beta_{0,q^F}=-\frac{1}{n}\left(\beta_{0,\chi_n,q}
-\frac{\chi_n(p)}{[p]_q}\beta_{0,\chi_n, q^p}\right).
\endaligned$$
Thus, we have $$\frac{\partial^n}{\partial
t^n}L_{p,q}(1-n,t|\chi)=-\frac{n!}{n}(p^*\frac{\log q}{q-1})^n
\left(\beta_{0,\chi_n,
q}-\chi_n(p)[p]_q^{-1}\beta_{0,\chi_n,q^p}\right).$$ Since
$\beta_{0,\chi, q}=0$ if $\chi \neq 1$, this become

$$\frac{\partial^n}{\partial t^n}L_{p,q}(1-n,t|\chi)=\cases
-(n-1)!(p^*\frac{\log
q}{q-1})^n(1-\frac{1}{p})\frac{q^F-1}{[F]\log q},&\text{ if }
\chi =1,\\
0,&\text{ if } \chi \neq 1.\endcases$$
 In the case $n=1 ,$ we easily see that
 $$\frac{\partial}{\partial t}L_{p,q}(0,t|\chi)=\cases
-p^*\frac{\log q}{q-1}(1-\frac{1}{p})\frac{q^F-1}{[F]\log
q},&\text{ if }
\chi =1,\\
0,&\text{ if } \chi \neq 1.\endcases$$ The value of
$\frac{\partial}{\partial s}L_{p,q}(0,t|\chi)$ is the coefficient
of $s$ in the expansion of $L_{p,q}(s,t|\chi)$ at $s=0 .$ By using
Taylor expansion at $s=0$, we see that
$$\aligned
&\frac{1}{1-s}=1+s+\cdots, \\
&<a+p^*t>^{1-s}=<a+p^*t>\left(1-s\log_p <a+p^* t>+\cdots \right),\\
&\binom{1-s}{m}=\frac{(-1)^{m+1}}{m(m-1)}s +\cdots .
\endaligned$$
By employing these expansion, along with some algebraic
manipulation, we evaluate $\frac{\partial}{\partial s}L_{p,q}(0,
t|\chi).$ From the definition of $L_{p,q}(s,t|\chi)$, we note that
$$\aligned
&L_{p,q}(s,t|\chi)\\
&=\frac{1}{s-1}\frac{1}{[F]_q}\sum_{\Sb a=1
\\(a, p)=1 \endSb}^F\chi(a)<a+p^*t>^{1-s}\sum_{m=0}^{\infty}
\binom{1-s}{m}\beta_{m,q^F}q^{(a+p^*t)m}[\frac{F}{a+p^*t}]_{q^{a+p^*t}}^m.
\endaligned$$
Thus, we have
$$\aligned
&\frac{\partial}{\partial s}L_{p,q}(s,t|\chi)|_{s=0}
=\sum_{a=1}^F\chi_1(a)\big\{\left(\frac{[a+p^*t]_q}{[F]_q}\beta_{0,q^F}+\beta_{1,q^F}\right)\log_p
<a+p^*t>\\
&-\frac{[a+p^*t]_q}{[F]_q}\beta_{0,q^F}
+\sum_{m=2}^{\infty}\frac{(-1)^m}{m(m-1)}
\frac{[F]_q^m}{[a+p^*t]_q^m}\beta_{m,q^F}\frac{[a+p^*t]_q}{[F]_q}\big\}\\
&+(q-1)\sum_{a=1}^F\chi_1(a)\big\{[a+p^*t]_q\beta_{1,q^F}\left(-1+\log_p
<a+p^*t>\right)\\
&+\sum_{m=2}^{\infty}\sum_{l=1}^m\binom ml
(q-1)^{l-1}[a+p^*t]_q^{l-m+1}[F]_q^{m-1}\beta_{m,q^F} \big\}.
\endaligned\tag26$$
We now define the Daehee $q$-operator, $D_{q,F}(x, y),$ as
follows:
$$D_{q,F}(x,y)=\left(\log_p x-1
\right)x\beta_{1,q}+\sum_{m=0}^{\infty}\sum_{l=1}^m\binom ml
\left(y-1\right)^{l-m}\left(y^F-1\right)^{m-1}x^{l-m+1}\beta_{m,q}
.\tag 27$$ In [7,8 ] the Diamond gamma function is defined by
$$G_p(x)=(x-\frac{1}{2})\log_p
x-x+\sum_{j=2}^{\infty}\frac{B_j}{j(j-1)}x^{1-j}, \text{ for
$|x|_p>1$ }.\tag28$$
 We now consider a $q$-analogue of the above Diamond gamma
 function as follows:
 $$G_{p,q}(x)=\int_{\Bbb
 Z_p}\left\{(x+[z]_q)\log_p(x+[z]_q)-(x+[z]_q)\right\}q^{-z}d\mu_q(z),
 \text{ for $|x|_p>1. $}\tag29 $$
From the above Eq.(29), we note that $G_{p,q}(x)$ is locally
analytic on $\Bbb C_p\setminus\Bbb Z_p $ . By (29), we easily see
that
$$G_{p,q}(x)=\left(x\beta_{0,q}+\beta_{1,q}\right)\log_p
x-x\beta_{0,q}+\sum_{n=1}^{\infty}\frac{(-1)^{n+1}}{n(n+1)}\beta_{n+1,q}x^{-n},\text{
for $|x|_p>1 $.}\tag30$$ Note that $\lim_{q\rightarrow
1}G_{p,q}(x)=G_p(x).$ Since, $w(a)$ is a root of unity for
$(a,p)=1$, we see  that
$$ \log_p<a+p^*t>=\log_p(a+p^*t)+\log_p w^{-1}(a)=\log_p(a+p^*t).
\tag 31$$ From the Eq.(26), Eq.(27), Eq.(30) and Eq.(31), we note
that
$$\aligned
&\frac{\partial}{\partial s}L_{p,q}(0,t|\chi) =\sum_{\Sb
a=1\\(a,p)=1\endSb}^F\chi_1(a)\big\{[F]_q^{-1}q^{p^*t}\log_p
[F]_q\beta_{0,q^F} [a]
+G_{p,q^F}\left(\frac{[a+p^*t]_q}{[F]_q}\right)\\
&+(q-1)D_{q^F,F}\left([a+p^*t]_q,q
\right)\big\}=-q^{p^*t}L_{p,q}(0,\chi)\log_p[F]\\
 &+\sum_{\Sb a=1
\\(a,p)=1\endSb}^F\chi_1(a)G_{p,q^F}\left(\frac{[a+p^*t]_q}{[F]_q}\right)
+(q-1) \sum_{\Sb a=1 \\(a,p)=1 \endSb}^F\chi_1(a)
D_{q^F,F}([a+p^*t]_q,q ).
\endaligned$$
Therefore we obtain the following theorem:
 \proclaim{ Theorem 7}
 Let $\chi$ be the primitive Dirichlet character, and let $F$ be a
 positive integral multiple of $p^*$ and $f=f_{\chi}$. Then for
 any $t\in\Bbb C_p$ with $|t|_p\leq 1$, we have
$$\aligned
\frac{\partial}{\partial s}L_{p,q}(0,t|\chi)& =\sum_{\Sb
a=1\\(a.p)=1\endSb}^F\chi_1(a)G_{p,q^F}\left(\frac{[a+p^*t]_q}{[F]_q}\right)
-q^{p^*t}L_{p,q}(0,\chi)\log_p [F]_q\\
&+(q-1)\sum_{\Sb a=1 \\ (a,p)=1 \endSb}^F\chi_1(a) D_{q^F,
F}\left([a+p^*t]_q, q \right).
\endaligned$$
 \endproclaim
Now we give the value of $L_{p,q}(s,t|\chi)$ at $s=1$ when $\chi
\neq 1$. From the definition of $L_{p,q}(s,t|\chi)$, we have
$$\aligned
L_{p,q}(s,t|\chi) &=\frac{1}{s-1}\frac{1}{[F]_q}\sum_{\Sb a=1\\
(a,p)=1
\endSb}^F\chi(a)\big\{<a+p^*t>^{1-s}\beta_{0,q^F}\\
&+<a+p^*t>^{1-s} \sum_{m=1}^{\infty}\binom{1-s}m
\beta_{m,q^F}q^{(a+p^*t)m}
\left[\frac{F}{a+p^*t}\right]_{q^{a+p^*t}}^m \big\}.
\endaligned$$
By using Taylor expansion  at $s=1$, we see that
$$\aligned
&\lim_{s\rightarrow 1}L_{p,q}(s,t|\chi)=-\frac{1}{[F]_q}\sum_{\Sb
a=1 \\(a,p)=1\endSb}^F\chi(a)\log_p <a+p^*t> \\
&-\frac{1}{[F]_q}\sum_{\Sb a=1
\\(a,p)=1\endSb}^F\chi(a)\sum_{m=1}^{\infty}(-1)^{m-1}\frac{(m-1)!}{m!}\beta_{m,q^F}q^{(a+p^*t)m}
\left(\frac{[F]_q}{[a+p^*t]_q}\right)^m .
\endaligned$$
Therefore we obtain the following theorem:
\proclaim{ Theorem 8}
Let $\chi$ be the Dirichlet character with conductor $f=f_{\chi}$
and let $F$ be the positive integral multiple of $p^*$ and
$f=f_{\chi}$. Then we have
$$\aligned
L_{p,q}(1,t|\chi) &=\frac{1}{[F]_q}\sum_{\Sb a=1 \\ (a,p)=1
\endSb}^F\chi(a)\big\{-\log_p
<a+p^*t>\\
&+\sum_{m=1}^{\infty}\frac{(-1)^{m-1}}{m} \beta_{m,
q^F}q^{(a+p^*t)m}\left(\frac{[f]_q}{[a+p^*t]_q}\right)^m\big\},
\text{ for $t\in\Bbb C_p$ with $|t|_p\leq 1$}.\endaligned$$
\endproclaim

Remark. From the above Theorem 8, we note that
$$\aligned
&L_{p,q}(1,0|\chi)=L_{p,q}(1,\chi)\\
&=\frac{1}{[F]_q}\sum_{\Sb a=1
\\(a,p)=1\endSb}\chi(a)\left(-\log_p
<a>+\sum_{m=1}^{\infty}\frac{(-1)^m}{m}\beta_{m,q^F}q^{am}\left(\frac{[F]_q}{[a]_q}\right)^m\right),
\endaligned$$
 and $\lim_{q\rightarrow 1} L_{p,q}(1,\chi)=L_p(1, \chi).$

\enddocument